\documentclass{amsart}
\title{Quasipositivity and braid index of pretzel knots}
\author{Lukas Lewark}
\address{\rule{0pt}{1cm}Faculty of Mathematics, University of Regensburg, 93053 Regensburg, Germany}
\email{\myemail{lukas@lewark.de}}
\urladdr{\url{http://www.lewark.de/lukas/}}
\keywords{Pretzel knots, quasipositive knots, Homflypt polynomial, braid index, Jones' conjecture, Khovanov-Rozansky homologies, slice-torus invariants}
\subjclass{57K10, 57K14, 57K18}
\usepackage{thmtools,enumerate,graphicx,xcolor,enumitem,afterpage,float,ucs,amsmath,amssymb,amscd,url,mathtools,microtype,tikz-cd,caption,amsthm,wrapfig}
\usepackage[latin1]{inputenc}
\usepackage{hyperref}
\usepackage[nameinlink]{cleveref}
\let\cref\Cref
\crefname{subsection}{subsection}{subsections}
\Crefname{subsection}{Subsection}{Subsections}
\Crefname{enumi}{}{}
\crefname{equation}{}{}
\definecolor{darkblue}{RGB}{0,0,96}
\definecolor{gray}{RGB}{127,127,127}
\definecolor{darkred}{RGB}{160,0,0}
\definecolor{lightyellow}{RGB}{255,255,128}
\hypersetup{colorlinks={true},linkcolor={black},citecolor={black},filecolor={black},urlcolor={black},pdfauthor={\authors},pdftitle={\shorttitle}}
\newcommand{\myemail}[1]{\href{mailto:#1}{#1}}
\newcommand{\qua}{\hskip 0.4em \ignorespaces}
\def\arxiv#1{\relax\ifhmode\unskip\qua\fi
\href{http://arxiv.org/abs/#1}%
{\tt arXiv:\penalty -100\unskip#1}}

\def\MR#1{\relax\ifhmode\unskip\qua\fi
\href{http://www.ams.org/mathscinet-getitem?mr=#1}{\tt MR#1}}
\def\xox#1{\csname xx#1\endcsname}

\renewenvironment{thebibliography}[1]{
  \begin{oldthebibliography}{#1}\small
    \setlength{\itemsep}{.5ex}
    \setlength{\parskip}{0em}
}
{
  \end{oldthebibliography}
}

\newcommand{\myqed}{\pushQED{\qed}\qedhere}

\declaretheorem{lemma}
\newtheorem{theorem}[lemma]{Theorem}

\newtheorem*{prize*}{Prize}
\newtheorem*{theorem*}{Theorem}

\newtheorem*{question*}{Question}
\theoremstyle{definition}
\newtheorem{definition}[lemma]{Definition}

\DeclareMathAlphabet{\mathpzc}{OT1}{pzc}{m}{it}
\DeclareMathOperator{\sgn}{sgn}
\DeclareMathOperator{\br}{br}

\newcommand{\Z}{\mathbb{Z}}

\graphicspath{{figures/}}
\begin{document}
\thispagestyle{empty}
\begin{abstract}
This short note is about three-stranded pretzel knots that have an even number of crossings in one of the strands.
We calculate the braid index of such knots and determine which of them are quasipositive.
The main tools are the Morton-Franks-Williams inequalities, and Khovanov-Rozansky concordance homomorphisms.

\end{abstract}
\maketitle
\section{Introduction}
\begin{wrapfigure}[10]{r}{58mm}
\centering
\raisebox{1mm}[0mm]{\includegraphics[angle=270,scale=.9]{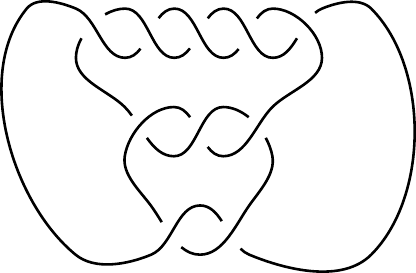}}
\captionsetup{width=.9\linewidth}
\caption{The $P(5,-3,-2)$ pretzel knot, which is quasipositive, but not strongly quasipositive.}
\label{fig:firstex}
\end{wrapfigure}

Our protagonists are the $P(p,q,-2r)$ pretzel knots, where $p,q,r$ are integers,
$p,q$ are odd and not $\pm 1$, and $r$ is not $0$. See \cref{fig:firstex} for an example.
Recently, Boileau, Boyer and Gordon proved that $P(p,q,-2r)$ is a strongly quasipositive knot
if and only if all of $p,q,r$ are positive \cite{MR4016557}.
The first result of this note gives a similar condition for these pretzel knots to be quasipositive.
\begin{theorem}
\label{prop:pretzqp}
Let $p,q$ be odd integers not equal to $\pm 1$, and $r$ a non-zero integer.
Then the $P(p,q,-2r)$ pretzel knot is quasipositive if and only if
$p + q\geq 0$ and $r > 0$.
\end{theorem}
As a second result, which we need to prove the first,
we explicitly define a 
braid $\beta(p,q,-2r)$ on $|r| + 2$ strands with closure
$P(p,q,-2r)$, and show that it is \emph{minimal},
i.e.~that it has the minimal number of strands among all braids with that closure.
\begin{theorem}
\label{prop:pretzbr}
Let $p,q$ be odd integers not equal to $\pm 1$, and $r$ a non-zero integer.
Then the braid index of $P(p,q,-2r)$ is $|r| + 2$. It is realized by the braid $\beta(p,q,-2r)$.
\end{theorem}
The braid $\beta(p,q,-2r)$ is defined in \cref{sec:br}.
To show its minimality and thus prove \cref{prop:pretzbr},
we partially compute the Homflypt polynomial of $P(p,q,-2r)$, and then rely on the Morton-Franks-Williams inequalities.

To prove the `if' direction of \cref{prop:pretzqp}, we simply observe that $\beta(p,q,-2r)$ is quasipositive
if $p+q\geq 0$ and $r>0$. To show the `only if' direction, we require the following obstruction to quasipositivity.
\begin{restatable}{lemma}{qpobs}\cite[Lemma~3.6]{MR3694648}
\label{lemma:qpobs}
Let $K$ be a quasipositive knot with braid index $b$. Let $w$ be the writhe of any minimal braid of $K$.
Let $\phi$ be any slice-torus invariant.  Then
\[
1 + w - b = 2\phi(K).\myqed
\]
\end{restatable}
The proof of \cref{lemma:qpobs} uses Jones' conjecture, shown in \cite{MR3235791,MR3302972}, stating that all minimal braids of a knot $K$ have the same writhe.
Comparing \cref{lemma:qpobs} with the previously computed values of the Khovanov-Rozansky $\mathfrak{sl}_3$-slice torus invariant~\cite{lew2} reveals that $P(p,q,-2r)$ is not quasipositive if $p+q<0$ or $r<0$.
The proof of \cref{prop:pretzqp} is contained in \cref{sec:qp}.

The canonical diagram $D(p,q,-2r)$ of the $P(p,q,-2r)$ pretzel consists of
two disks connected by three bands with $p,q,-2r$ half-twists, respectively.
Note that by convention the signs of $p,q,-2r$ encode the \emph{handedness} of twists.
The \emph{signs} of the crossings in the three bands in $D(p,q,-2r)$
are $\sgn(p), \sgn(q), \sgn(r)$, respectively.
As further background, let us list some results in a similar vein as \cref{prop:pretzqp}:\medskip
\begin{itemize}
\setlength\itemsep{.5ex}
\item $P(p,q,-2r)$ is strongly quasipositive $\Leftrightarrow$ $P(p,q,-2r)$ is positive\\$\Leftrightarrow$ $D(p,q,-2r)$ is positive $\Leftrightarrow$ $p,q,r$ are all positive \cite{MR4016557}.
\item $P(p,q,-2r)$ is alternating $\Leftrightarrow$ $D(p,q,-2r)$ is alternating\\$\Leftrightarrow$ $p,q,-r$ all have the same sign \cite{MR966948}.
\item $P(p,q,-2r)$ is quasi-alternating $\Leftrightarrow$
$p,q,-r$ all have the same sign, or\break $\{p+q,p-2r,q-2r\}$ contains a negative and a positive number \cite{greene2}.
\item $P(p,q,-2r)$ is fibred $\Leftrightarrow$ $p,q$ are of opposite sign, or $|r| = 1$ \cite{gabai2}.
\item Conjecturally, $P(p,q,-2r)$ is slice $\Leftrightarrow$ $p + q = 0$ \cite{MR3402337}.
\end{itemize}\medskip

\paragraph{\bfseries Acknowledgments.}
I warmly thank Peter Feller and Andrew Lobb for the inspiring collaboration on squeezed knots \cite{sqz} and slice-torus invariants \cite{sti},
from which this note originated.
I gratefully acknowledge support by the DFG, project no.~412851057.
Thanks to Paula Tru\"ol for comments on a first version of this note.

\section{The braid index of pretzel knots}
\label{sec:br}
\begin{figure}[t]
\begin{center}
\includegraphics[angle=-90,scale=.9]{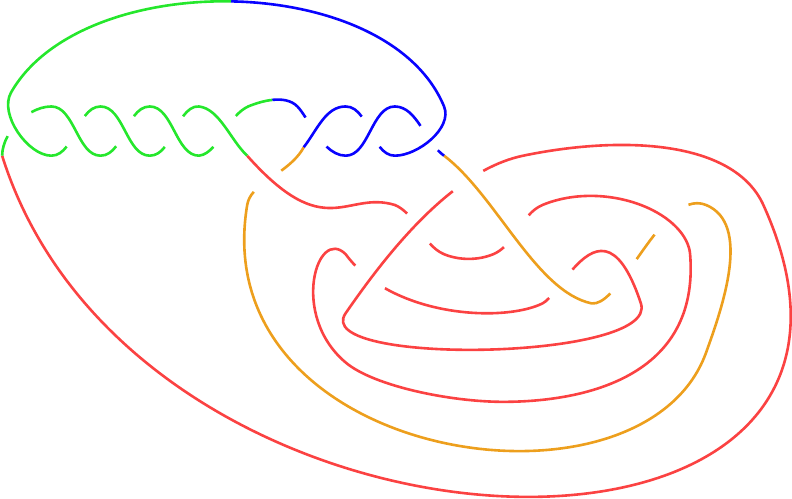}\qquad
\includegraphics[angle=-90,scale=.9]{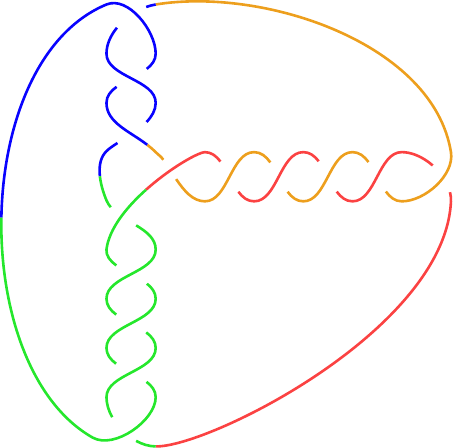}
\end{center}
\caption{On the left, the closure of $\beta(5,-3,-6) = a_1^5 a_2 a_1^{-3} a_4^{-1}a_3^{-1} a_2 a_3 a_4^2 a_3$. On the right, the standard diagram $D(5,-3,-6)$ of the pretzel knot $P(5,-3,-6)$.
The reader is invited to spot the isotopy between the two diagrams.}
\label{fig:pretzelbraid}
\end{figure}

\begin{definition}\label{def:beta}
Let $p,q$ be odd integers not equal to $\pm 1$, and $r$ a non-zero integer.
Let us define the following braids on $|r| + 2$ strands,
denoting the standard generators of the braid group by $a_1, \ldots, a_{|r|+1}$:
\begin{align*}
\gamma_r & = a_3 a_4 \cdots a_{|r|} a_{|r| + 1},\\
\overline{\gamma_r} & = a_{|r| + 1} a_{|r|} \cdots a_4 a_3, \\
\beta(p,q,-2r) & = a_1^p a_2 a_1^q \gamma_r^{-1} a_2 \gamma_r \overline{\gamma_r} \quad\text{if $r > 0$,}\\
\beta(p,q,-2r) & = a_1^p a_2^{-1} a_1^q \overline{\gamma_r} a_2^{-1}\overline{\gamma_r}^{-1} \gamma_r^{-1} \quad\text{if $r < 0$.}
\end{align*}
\end{definition}
It has been shown in \cite{brpr} that the closure of $\beta(p,q,-2r)$ is $P(p,q,-2r)$ (see also \Cref{fig:pretzelbraid} for an example).
Note that in the simplest cases $r = \pm 1$, $\gamma_r$ and $\overline{\gamma_r}$ are empty and thus $\beta(p,q,-2) = a_1^p a_2 a_1^q a_2$ and $\beta(p,q,2) = a_1^p a_2^{-1} a_1^q a_2^{-1}$.

Let us now prove \cref{prop:pretzbr} by showing that $\beta(p,q,-2r)$ realizes the braid index; for this, we will use the Morton-Franks-Williams inequalities \cite{MR809504,MR896009}:
\begin{equation}\label{eq:mfw1}
1 + w(\beta) - k(\beta) \leq e(K) \leq E(K) \leq -1 + w(\beta) + k(\beta),
\end{equation}
where $\beta$ is a braid with writhe $w(\beta)$, on $k(\beta)$ strands, and closure a knot $K$.
Moreover, $e$ and $E$ are respectively the minimum and maximum exponent of $v$ appearing in the Homflypt polynomial
$Q_K(v,z) \in \Z[v^{\pm 1}, z^{\pm 1}]$, which is defined by the skein relation
$v^{-1} Q_{L_+} - v Q_{L_-} = zQ_{L_0}$ 
 and setting $Q_U = 1$ (where $L_+$ is a link admitting a diagram $D$, such that $L_-$ arises from $D$ by changing a positive crossing to a negative one, and $L_0$ by orientedly resolving that crossing).
Now \eqref{eq:mfw1} implies
\[2\br(K) \geq E(K) - e(K) + 2,\]
where $\br(K)$ denotes the braid index of $K$.
For our purposes, it is sufficient to use the specialization
$Q'_K(v) \coloneqq Q_K(v,1) \in \Z[v^{\pm 1}]$ of the Homflypt polynomial. Defining $e'$ and $E'$
as minimum and maximum exponent of $v$ appearing in $Q'$, we clearly have $E'\leq E$ and $e'\geq e$, and so
\begin{equation} \label{eq:mfw2}
2\br(K) \geq E'(K) - e'(K) + 2.
\end{equation}
This is the lower bound for the braid index that we will use. Let us now compute $e'$ and $E'$ of
$P(p,q,-2r)$.
We start by computing $Q'$ of the torus link $T(2,n)$ for all~$n\in \Z$,
which is the closure of the two-stranded braid $a_1^n$. One finds
\begin{align*}
Q'_{T(2,0)} & = v^{-1} - v, \\
Q'_{T(2,1)} & = 1, \\
Q'_{T(2,n+2)} & = v Q'_{T(2,n+1)} + v^2Q'_{T(2,n)} \qquad\text{for all $n \geq 0$}.
\end{align*}
Denote by $F(n)$ the $n$-th Fibonacci number, i.e.\ $F(0) = 0, F(1) = 1$ and $F(n+2) = F(n+1) + F(n)$ for all $n\geq 0$.
It easily follows inductively that for all $n \geq 1$,
\begin{equation}\label{eq:homflypttorus}
Q'_{T(2,n)} = F(n+1) v^{n-1} - F(n - 1) v^{n+1}.
\end{equation}
Using $Q'_{T(2,-n)}(v) = (-1)^{n+1}Q'_{T(2,n)}(v^{-1})$, one may extend \eqref{eq:homflypttorus}
to all $n$ by setting $F(-n) = (-1)^{n+1} F(n)$ for all positive $n$. Since $F(n) \neq 0$ for all $n\neq 0$,
we find $e'(T(2,n)) = n-1$ and $E'(T(2,n))=n+1$ for all $n\neq \pm 1$.

Let us proceed to the pretzels, starting with the case $r = 1$.
Applying the skein relation to one of the two crossings in the band with 2 twists in $D(p,q,-2)$,
and using that $Q'(K\# J) = Q'(K)Q'(J)$ for knots $K,J$ yields
\[
Q'_{P(p,q,-2)} = v Q'_{T(2,p+q)} + v^2 Q'_{T(2,p)} Q'_{T(2,q)}.
\]
In this polynomial, the term with highest exponent of $v$ is $F(p-1)F(q-1)v^{p+q+4}$.
The term with lowest exponent of $v$ is
$x v^{p+q}$, provided
\begin{equation}\label{eq:fibonaccidetail}
x \coloneqq F(p+q+1) + F(p+1)F(q+1) \neq 0.
\end{equation}
Let us show \cref{eq:fibonaccidetail}. We have $F(p+q+1) > 0$,
$\sgn F(p+1) = \sgn p, \sgn F(q+1) = \sgn q$. So if $\sgn p = \sgn q$, then $F(p+1)F(q+1) > 0$, and \cref{eq:fibonaccidetail} follows.
If $\sgn p \neq \sgn q$, then $|p + 1| > |p + q + 1|$ or $|q + 1| > |p + q + 1|$, and 
thus $|F(p+1)| > |F(p + q + 1)|$ or $|F(q+1)| > |F(p + q + 1)|$ (using monotonicity of the Fibonacci numbers), and  so \eqref{eq:fibonaccidetail} also holds.
It follows that $e'(P(p,q,-2)) = p+q$ and $E'(P(p,q,-2)) = p+q+4$, so $\br(P(p,q,-2)) \geq 3 = r + 2$ as desired.

Now for the case $r \geq 2$, applying the skein relation to one of the crossings in the band with $2r$ twists in $D(p,q,-2r)$ gives
\[
Q'_{P(p,q,-2r)} = v Q'_{T(2,p+q)} + v^{2} Q'_{P(p,q,2-2r)}.
\]
It follows inductively that $Q'_{P(p,q,-2r)} = A + B$ for
\begin{align*}
 A &= (v + v^{3} + \ldots + v^{2r-1})\cdot Q'_{T(2,p+q)},\\
 B &= v^{2r}\cdot Q'_{T(2,p)}\cdot Q'_{T(2,q)}
\end{align*}
and so
\begin{alignat*}{3}
e'(A) & = p + q, & E'(A) & = p+q+2r \\
e'(B) & = p+q+2r-2,\quad & E'(B) & = p+q+2r+2.
\end{alignat*}
Thus $e'(P(p,q,-2r)) = e'(A)$ and $E'(P(p,q,-2r)) = E'(B)$, and it follows
that $\br(P(p,q,-2r)) \geq r + 2$.

We have shown for all positive $r$ that the braid index of $P(p,q,-2r)$ equals $|r| + 2$.
The braid index of the knots $P(p,q,-2r)$ and $P(-p,-q,2r)$ agrees, since they are mirror images of one another.
Thus we have completed the proof of \cref{prop:pretzbr}.
\myqed

\section{The quasipositivity of pretzel knots}
\label{sec:qp}
As sketched in the introduction, the proof of \cref{prop:pretzqp} relies on \cref{lemma:qpobs},
which we restate below for the reader's convenience.
Here, a braid is \emph{quasipositive} if it equals a product of conjugates of the standard generators $a_i$ of the braid group. A knot is \emph{quasipositive} it it is the closure of some quasipositive braid.
A \emph{slice-torus invariant} $\phi$ is a homomorphism from the smooth knot concordance group to $\mathbb{R}$,
such that $\phi(K) \leq g_4(K)$ holds for all knots $K$ (where $g_4$ denotes the smooth slice genus),
and $\phi(T(p,q)) = g_4(T(p,q))$ holds for all positive torus knots $T(p,q)$ \cite{livingston} (see also \cite{lew2,sti}).

\qpobs*
While the proof of \cref{lemma:qpobs} uses Jones' conjecture, 
it is worth noting that if the Morton-Franks-Williams inequality \cref{eq:mfw2} is an equality for a knot $K$
(as it is for the pretzel knots we are considering), then the statement of Jones' conjecture for $K$
can be easily deduced directly from \cref{eq:mfw1}.

Let us now prove \cref{prop:pretzqp}. We consider four exhaustive and mutually exclusive cases.
Note that the writhe $w(\beta(p,q,-2r))$ equals $p + q + r + \sgn r$.
We will use the statement of \cref{prop:pretzbr} that $\beta(p,q,-2r)$ is a minimal braid representative for~$P(p,q,-2r)$.
\begin{itemize}
 \setlength\itemsep{1ex}
\item Case $p+q\geq 0$ and $r>0$.
Then, $\beta(p,q,-2r)$ is clearly quasipositive.
This settles the `if' direction of \cref{prop:pretzqp}. The remaining three cases cover the `only if' direction.
\item Case $p+q\leq 0$ and $r<0$.
Since $w(\beta(p,q,-2r))$ is negative, it follows that $1 + w(\beta(p,q,-2r)) - \br(P(p,q,-2r))$ is also negative. So by \cref{lemma:qpobs}, $\phi(P(p,q,-2r)) < 0$. However, quasipositive knots have non-negative slice-torus invariants \cite{lew2},
so $P(p,q,-2r)$ is not quasipositive.
Note that $P(p,q,-2r)$ is in fact quasinegative; so this case also follows from Hayden's theorem \cite{MR3788795} that
no non-trivial knot is both quasipositive and quasinegative.
\item Case $p + q < 0$ and $r > 0$. 
If $P(p,q,-2r)$ is quasipositive, by \cref{lemma:qpobs}, it follows that any slice-torus invariant $\phi$ satisfies
$2\phi(P(p,q,-2r)) = p+q$.  However, as computed in \cite{lew2}, the Khovanov-Rozansky $\mathfrak{sl}_3$-slice torus invariant $s_3$ (which takes values in $\tfrac12\mathbb{Z}$) satisfies
\[
2s_3(P(p,q,-2r)) \in \{p+q+2, p + q+1\},
\]
contradicting the assumption of quasipositivity of $P(p,q,-2r)$.
\item Case $p + q > 0$ and $r < 0$.
If $P(p,q,-2r)$ were quasipositive, then by \cref{lemma:qpobs} any slice torus invariant $\phi$ would satisfy
$2\phi(K) = p + q + 2r - 2$.
However,
\[
2s_3(P(p,q,-2r)) \in \{p+q-2, p + q-1\},
\]
leading to a contradiction.\myqed
\end{itemize}

\section{Further properties of pretzel knots}
Let us have a quick look at further notions of positivity. A knot is called \emph{braid positive} if it is the closure
of a positive braid, i.e.~a braid that can be written as product of positive powers of the standard generators.
Such knots are fibered~\cite{MR520522} and positive.
Recalling from the introduction which $P(p,q,-2r)$ pretzel knots are fibered, and which are positive,
one sees that for $P(p,q,-2r)$ to be braid positive, it is necessary that $p$ and $q$ are positive and $r = 1$.
Since $\beta(p,q,-2) = a_2^p a_1 a_2^q a_1$, this condition is also sufficient.

A knot is called \emph{strongly quasipositive} if it is the closure of a \emph{strongly quasipositive braid},
i.e.~a braid on $n$ strands that can be written as product of certain conjugates of the standard generators, namely
$b_{ij} \coloneqq (a_{i+1} \cdots a_j)^{-1} a_i (a_{i+1} \cdots a_j)$ for all\break $1 \leq i \leq j < n$.
Note that $b_{ii} = a_i$.
Now, recall that $P(p,q,-2r)$ is strongly quasipositive if and only if $p, q, r > 0$.
On the other hand, one observes that for $p,q,r>0$, the braid $\beta(p,q,-2r) = a_1^p a_2 a_1^q b_{2,r+1} \overline{\gamma_r}$ is strongly quasipositive.

In summary we have:\medskip
\begin{itemize}
\setlength\itemsep{.5ex}
\item $P(p,q,-2r)$ braid positive $\Leftrightarrow$ $\beta(p,q,-2r)$ positive braid $\Leftrightarrow$ $p, q > 0$, $r = 1$. 
\item $P(p,q,-2r)$ str.~quasipos.~$\Leftrightarrow$ $\beta(p,q,-2r)$ str.~quasipos.~$\Leftrightarrow$ $p, q, r >0$.
\item $P(p,q,-2r)$ quasipositive $\Leftrightarrow$ $\beta(p,q,-2r)$ quasipositive $\Leftrightarrow$ $p + q\geq 0$, $r > 0$.
\end{itemize}\medskip

A knot $K$ is called \emph{squeezed} if it appears as a slice of a genus-minimizing oriented connected smooth cobordism between a positive and a negative torus knot \cite{sqz}. The class of squeezed knots includes quasipositive, quasinegative, and alternating knots.
This guarantess the squeezedness of $P(p,q,-2r)$ unless $\sgn (p + q) = - \sgn r$ and $\sgn p = -\sgn q$.
For pretzels satisfying those equations, let us distinguish two cases.
In the first case that $\sgn (p - 2r) = - \sgn (q - 2r)$,
e.g.~as for $P(5,-3,2)$,
the Khovanov-Rozansky $\mathfrak{sl}_3$-slice torus invariant $s_3$
is not equal to the Rasmussen invariant \cite{lew2}. Since 
all slice-torus invariants take the same value on a fixed squeezed knot, the knots in that case are not squeezed.
I have not been able to determine the squeezedness of the pretzel knots
of the second case, namely those with $\sgn (p - 2r) = \sgn (q - 2r)$, such as~$P(5,-3,4)$.

\bibliographystyle{myamsalpha}
\bibliography{References}
\end{document}